\documentclass{gtart}

\def\ifplaintex{\expandafter\ifx\csname documentclass\endcsname\relax}

\def\gtp{{\mathsurround=0pt\it $\cal G\mskip-2mu$eometry \&\ 
$\cal T\!\!$opology $\cal P\!$ublications}}  

\def\recd{{\small Received:\qua\receiveddate\ifx\reviseddate\relax
\else\qquad Revised:\qua\reviseddate\fi\par}} 

\def\lognumber#1{\def\thelognumber{#1}}
\def\volumenumber#1{\def\thevolumenumber{#1}}
\def\volumeyear#1{\def\thevolumeyear{#1}}
\def\papernumber#1{\def\thepapernumber{#1}}
\def\pagenumbers#1#2{\def\startpage{#1}\def\finishpage{#2}}
\def\published#1{\def\publishdate{#1}}

\def\received#1{\def\receiveddate{#1}}
\def\revised#1{\def\reviseddate{#1}}
\def\accepted#1{\def\accepteddate{#1}}

\long\def\asciiabstract#1{\long\def\theasciiabstract{#1}}
\def\asciikeywords#1{\def\theasciikeywords{#1}}

\let\\\par\let\thelognumber\relax\let\thevolumenumber\relax
\let\thepapernumber\relax\let\thevolumeyear\relax\let\startpage\relax
\let\finishpage\relax\let\publishdate\relax\let\receiveddate\relax
\let\reviseddate\relax\let\accepteddate\relax\let\theasciititle\relax
\let\theasciiauthors\relax
\let\theasciiabstract\relax\let\theasciikeywords\relax

\let\theasciiemail\relax

\ifplaintex
\font\logobig=cmssbx10 scaled 3836
\font\logomed=cmssbx10 scaled 2557
\else
\font\logobig=cmssbx10 scaled 4200
\font\logomed=cmssbx10 scaled 2800
\fi

\long\def\makeagttitle{   
\count0=\startpage
\agt\hfill      
\hbox to 45truept{\vbox to 0pt{\vglue -13truept{\logomed A\kern -.37em{\logobig 
T}\kern -.38em G}\vss}\hss}
\break
{\small Volume \thevolumenumber\ (\thevolumeyear)
\startpage--\finishpage\nl
Published: \publishdate}

\vglue .25truein

{\parskip=0pt\leftskip 0pt plus
1fil\def\\{\par\smallskip}{\Large\bf\thetitle}\par\medskip} \vglue
0.05truein

{\parskip=0pt\leftskip 0pt plus 1fil\def\\{\par}{\sc\theauthors}
\par\medskip}%
 
\vglue 0.03truein 

{\small\leftskip 25truept\rightskip 25truept{\bf Abstract}\stdspace\theabstract

{\bf AMS Classification}\stdspace\theprimaryclass
\ifx\thesecondaryclass\relax\else; \thesecondaryclass\fi\par
{\bf Keywords}\stdspace \thekeywords\par}\vglue 7truept

}   

\ifplaintex
\font\phead=cmsl9 scaled 950
\font\pnum=cmbx10 scaled 913
\font\pfoot=cmsl9 scaled 950
\headline{\vbox to 0pt{\vskip -4.5mm\line{\small\phead\ifnum
\count0=\startpage ISSN 1472-2739 (on-line) 1472-2747 (printed)
\hfill {\pnum\folio}\else\ifodd\count0\def\\{ }%
\ifx\theshorttitle\relax\thetitle\else\theshorttitle\fi\hfill{\pnum\folio}
\else\def\\{ and }{\pnum\folio}\hfill\ifx\theshortauthors\relax\theauthors
\else\theshortauthors\fi\fi\fi}\vss}}
\footline{\vbox to 0pt{\vglue 0mm\line{\small\pfoot\ifnum\count0=\startpage
\copyright\ \gtp\hfill\else
\agt, Volume \thevolumenumber\ (\thevolumeyear)\hfill\fi}\vss}}
\else
\font=cmsl9 scaled 1050
\font\lnum=cmbx10 
\font=cmsl9 scaled 1050
\makeatletter
\def\@oddhead{{\small\lhead\ifnum\count0=\startpage ISSN 1472-2739 
(on-line) 1472-2747 (printed)\hfill {\lnum\number\count0}\else\ifodd\count0
\def\\{ }\ifx\theshorttitle\relax \thetitle \else\theshorttitle\fi\hfill
{\lnum\number\count0}\else\def\\{ and }{\lnum\number\count0}
\hfill\ifx\theshortauthors\relax 
\theauthors\else\theshortauthors\fi\fi\fi}}\def\@evenhead{\@oddhead}
\def\@oddfoot{\small\lfoot\ifnum\count0=\startpage\copyright\ \gtp\hfill\else
\agt, Volume \thevolumenumber\ (\thevolumeyear)\hfill\fi}
\def\@evenfoot{\@oddfoot}
\makeatother
\fi
\let\maketitlepage\makeagttitle
\let\makeshorttitle\maketitlepage
\let\maketitle\maketitlepage

\newwrite\gtoutfile
\long\gdef\makeheadfile{  
{\def\\{, }\def\s{ }
\immediate\openout\gtoutfile head.xxx
\immediate\write\gtoutfile{To: math@arxiv.org}
\immediate\write\gtoutfile{Subject: put OR rep NNNNN:ppppp}
\immediate\write\gtoutfile{--text follows this line--}
\immediate\write\gtoutfile{Proxy-for: \ifx\theasciiauthors\relax
\theauthors\else\theasciiauthors\fi\s<\ifx\theasciiemail\relax\theemail\else\theasciiemail\fi>}
\immediate\write\gtoutfile{\noexpand\\}
\immediate\write\gtoutfile{Authors: \ifx\theasciiauthors\relax
\theauthors\else\theasciiauthors\fi}
{\def\\{ }\immediate\write\gtoutfile{Title: \ifx\theasciititle\relax
\thetitle\else\theasciititle\fi}}
\immediate\write\gtoutfile{Subj-class: GT or SG, GR etc}
\immediate\write\gtoutfile{MSC-class: \theprimaryclass\ifx\thesecondaryclass\relax\else, \thesecondaryclass\fi}
\immediate\write\gtoutfile{Journal-ref: Algebr. Geom. Topol. \thevolumenumber\s
(\thevolumeyear) \startpage-\finishpage}
\immediate\write\gtoutfile{Comments: Published by Algebraic and
Geometric Topology at}
\immediate\write\gtoutfile{\s\s\s  http://www.maths.warwick.ac.uk/agt/AGTVol\thevolumenumber/agt-\thevolumenumber-\thepapernumber.abs.html}
\immediate\write\gtoutfile{\noexpand\\}
\immediate\write\gtoutfile{}
\ifx\theasciiabstract\relax
\immediate\write\gtoutfile{\theabstract}\else
\immediate\write\gtoutfile{\theasciiabstract}\fi
\immediate\write\gtoutfile{}
\immediate\write\gtoutfile{\noexpand\\}
\immediate\write\gtoutfile{}
\immediate\closeout\gtoutfile}}  

\def\maketitlepage{\makeagttitle\makeheadfile}
\let\makeshorttitle\maketitlepage
\let\maketitle\maketitlepage

\def\ifplaintex{\expandafter\ifx\csname documentclass\endcsname\relax}

\def\gtp{{\mathsurround=0pt\it $\cal G\mskip-2mu$eometry \&\ 
$\cal T\!\!$opology $\cal P\!$ublications}}  

\def\recd{{\small Received:\qua\receiveddate\ifx\reviseddate\relax
\else\qquad Revised:\qua\reviseddate\fi\par}} 

\def\lognumber#1{\def\thelognumber{#1}}
\def\volumenumber#1{\def\thevolumenumber{#1}}
\def\volumeyear#1{\def\thevolumeyear{#1}}
\def\papernumber#1{\def\thepapernumber{#1}}
\def\pagenumbers#1#2{\def\startpage{#1}\def\finishpage{#2}}
\def\published#1{\def\publishdate{#1}}

\def\received#1{\def\receiveddate{#1}}
\def\revised#1{\def\reviseddate{#1}}
\def\accepted#1{\def\accepteddate{#1}}

\long\def\asciiabstract#1{\long\def\theasciiabstract{#1}}
\def\asciikeywords#1{\def\theasciikeywords{#1}}

\let\\\par\let\thelognumber\relax\let\thevolumenumber\relax
\let\thepapernumber\relax\let\thevolumeyear\relax\let\startpage\relax
\let\finishpage\relax\let\publishdate\relax\let\receiveddate\relax
\let\reviseddate\relax\let\accepteddate\relax\let\theasciititle\relax
\let\theasciiauthors\relax
\let\theasciiabstract\relax\let\theasciikeywords\relax

\let\theasciiemail\relax

\ifplaintex
\font\logobig=cmssbx10 scaled 3836
\font\logomed=cmssbx10 scaled 2557
\else
\font\logobig=cmssbx10 scaled 4200
\font\logomed=cmssbx10 scaled 2800
\fi

\long\def\makeagttitle{   
\count0=\startpage
\agt\hfill      
\hbox to 45truept{\vbox to 0pt{\vglue -13truept{\logomed A\kern -.37em{\logobig 
T}\kern -.38em G}\vss}\hss}
\break
{\small Volume \thevolumenumber\ (\thevolumeyear)
\startpage--\finishpage\nl
Published: \publishdate}

\vglue .25truein

{\parskip=0pt\leftskip 0pt plus
1fil\def\\{\par\smallskip}{\Large\bf\thetitle}\par\medskip} \vglue
0.05truein

{\parskip=0pt\leftskip 0pt plus 1fil\def\\{\par}{\sc\theauthors}
\par\medskip}%
 
\vglue 0.03truein 

{\small\leftskip 25truept\rightskip 25truept{\bf Abstract}\stdspace\theabstract

{\bf AMS Classification}\stdspace\theprimaryclass
\ifx\thesecondaryclass\relax\else; \thesecondaryclass\fi\par
{\bf Keywords}\stdspace \thekeywords\par}\vglue 7truept

}   

\ifplaintex
\font\phead=cmsl9 scaled 950
\font\pnum=cmbx10 scaled 913
\font\pfoot=cmsl9 scaled 950
\headline{\vbox to 0pt{\vskip -4.5mm\line{\small\phead\ifnum
\count0=\startpage ISSN 1472-2739 (on-line) 1472-2747 (printed)
\hfill {\pnum\folio}\else\ifodd\count0\def\\{ }%
\ifx\theshorttitle\relax\thetitle\else\theshorttitle\fi\hfill{\pnum\folio}
\else\def\\{ and }{\pnum\folio}\hfill\ifx\theshortauthors\relax\theauthors
\else\theshortauthors\fi\fi\fi}\vss}}
\footline{\vbox to 0pt{\vglue 0mm\line{\small\pfoot\ifnum\count0=\startpage
\copyright\ \gtp\hfill\else
\agt, Volume \thevolumenumber\ (\thevolumeyear)\hfill\fi}\vss}}
\else
\font=cmsl9 scaled 1050
\font\lnum=cmbx10 
\font=cmsl9 scaled 1050
\makeatletter
\def\@oddhead{{\small\lhead\ifnum\count0=\startpage ISSN 1472-2739 
(on-line) 1472-2747 (printed)\hfill {\lnum\number\count0}\else\ifodd\count0
\def\\{ }\ifx\theshorttitle\relax \thetitle \else\theshorttitle\fi\hfill
{\lnum\number\count0}\else\def\\{ and }{\lnum\number\count0}
\hfill\ifx\theshortauthors\relax 
\theauthors\else\theshortauthors\fi\fi\fi}}\def\@evenhead{\@oddhead}
\def\@oddfoot{\small\lfoot\ifnum\count0=\startpage\copyright\ \gtp\hfill\else
\agt, Volume \thevolumenumber\ (\thevolumeyear)\hfill\fi}
\def\@evenfoot{\@oddfoot}
\makeatother
\fi
\let\maketitlepage\makeagttitle
\let\makeshorttitle\maketitlepage
\let\maketitle\maketitlepage

\newwrite\gtoutfile
\long\gdef\makeheadfile{  
{\def\\{, }\def\s{ }
\immediate\openout\gtoutfile head.xxx
\immediate\write\gtoutfile{To: math@arxiv.org}
\immediate\write\gtoutfile{Subject: put OR rep NNNNN:ppppp}
\immediate\write\gtoutfile{--text follows this line--}
\immediate\write\gtoutfile{Proxy-for: \ifx\theasciiauthors\relax
\theauthors\else\theasciiauthors\fi\s<\ifx\theasciiemail\relax\theemail\else\theasciiemail\fi>}
\immediate\write\gtoutfile{\noexpand\\}
\immediate\write\gtoutfile{Authors: \ifx\theasciiauthors\relax
\theauthors\else\theasciiauthors\fi}
{\def\\{ }\immediate\write\gtoutfile{Title: \ifx\theasciititle\relax
\thetitle\else\theasciititle\fi}}
\immediate\write\gtoutfile{Subj-class: GT or SG, GR etc}
\immediate\write\gtoutfile{MSC-class: \theprimaryclass\ifx\thesecondaryclass\relax\else, \thesecondaryclass\fi}
\immediate\write\gtoutfile{Journal-ref: Algebr. Geom. Topol. \thevolumenumber\s
(\thevolumeyear) \startpage-\finishpage}
\immediate\write\gtoutfile{Comments: Published by Algebraic and
Geometric Topology at}
\immediate\write\gtoutfile{\s\s\s  http://www.maths.warwick.ac.uk/agt/AGTVol\thevolumenumber/agt-\thevolumenumber-\thepapernumber.abs.html}
\immediate\write\gtoutfile{\noexpand\\}
\immediate\write\gtoutfile{}
\ifx\theasciiabstract\relax
\immediate\write\gtoutfile{\theabstract}\else
\immediate\write\gtoutfile{\theasciiabstract}\fi
\immediate\write\gtoutfile{}
\immediate\write\gtoutfile{\noexpand\\}
\immediate\write\gtoutfile{}
\immediate\closeout\gtoutfile}}  

\def\maketitlepage{\makeagttitle\makeheadfile}
\let\makeshorttitle\maketitlepage
\let\maketitle\maketitlepage

\def\ifplaintex{\expandafter\ifx\csname documentclass\endcsname\relax}

\def\gtp{{\mathsurround=0pt\it $\cal G\mskip-2mu$eometry \&\ 
$\cal T\!\!$opology $\cal P\!$ublications}}  

\def\recd{{\small Received:\qua\receiveddate\ifx\reviseddate\relax
\else\qquad Revised:\qua\reviseddate\fi\par}} 

\def\lognumber#1{\def\thelognumber{#1}}
\def\volumenumber#1{\def\thevolumenumber{#1}}
\def\volumeyear#1{\def\thevolumeyear{#1}}
\def\papernumber#1{\def\thepapernumber{#1}}
\def\pagenumbers#1#2{\def\startpage{#1}\def\finishpage{#2}}
\def\published#1{\def\publishdate{#1}}

\def\received#1{\def\receiveddate{#1}}
\def\revised#1{\def\reviseddate{#1}}
\def\accepted#1{\def\accepteddate{#1}}

\long\def\asciiabstract#1{\long\def\theasciiabstract{#1}}
\def\asciikeywords#1{\def\theasciikeywords{#1}}

\let\\\par\let\thelognumber\relax\let\thevolumenumber\relax
\let\thepapernumber\relax\let\thevolumeyear\relax\let\startpage\relax
\let\finishpage\relax\let\publishdate\relax\let\receiveddate\relax
\let\reviseddate\relax\let\accepteddate\relax\let\theasciititle\relax
\let\theasciiauthors\relax
\let\theasciiabstract\relax\let\theasciikeywords\relax

\let\theasciiemail\relax

\ifplaintex
\font\logobig=cmssbx10 scaled 3836
\font\logomed=cmssbx10 scaled 2557
\else
\font\logobig=cmssbx10 scaled 4200
\font\logomed=cmssbx10 scaled 2800
\fi

\long\def\makeagttitle{   
\count0=\startpage
\agt\hfill      
\hbox to 45truept{\vbox to 0pt{\vglue -13truept{\logomed A\kern -.37em{\logobig 
T}\kern -.38em G}\vss}\hss}
\break
{\small Volume \thevolumenumber\ (\thevolumeyear)
\startpage--\finishpage\nl
Published: \publishdate}

\vglue .25truein

{\parskip=0pt\leftskip 0pt plus
1fil\def\\{\par\smallskip}{\Large\bf\thetitle}\par\medskip} \vglue
0.05truein

{\parskip=0pt\leftskip 0pt plus 1fil\def\\{\par}{\sc\theauthors}
\par\medskip}%
 
\vglue 0.03truein 

{\small\leftskip 25truept\rightskip 25truept{\bf Abstract}\stdspace\theabstract

{\bf AMS Classification}\stdspace\theprimaryclass
\ifx\thesecondaryclass\relax\else; \thesecondaryclass\fi\par
{\bf Keywords}\stdspace \thekeywords\par}\vglue 7truept

}   

\ifplaintex
\font\phead=cmsl9 scaled 950
\font\pnum=cmbx10 scaled 913
\font\pfoot=cmsl9 scaled 950
\headline{\vbox to 0pt{\vskip -4.5mm\line{\small\phead\ifnum
\count0=\startpage ISSN 1472-2739 (on-line) 1472-2747 (printed)
\hfill {\pnum\folio}\else\ifodd\count0\def\\{ }%
\ifx\theshorttitle\relax\thetitle\else\theshorttitle\fi\hfill{\pnum\folio}
\else\def\\{ and }{\pnum\folio}\hfill\ifx\theshortauthors\relax\theauthors
\else\theshortauthors\fi\fi\fi}\vss}}
\footline{\vbox to 0pt{\vglue 0mm\line{\small\pfoot\ifnum\count0=\startpage
\copyright\ \gtp\hfill\else
\agt, Volume \thevolumenumber\ (\thevolumeyear)\hfill\fi}\vss}}
\else
\font=cmsl9 scaled 1050
\font\lnum=cmbx10 
\font=cmsl9 scaled 1050
\makeatletter
\def\@oddhead{{\small\lhead\ifnum\count0=\startpage ISSN 1472-2739 
(on-line) 1472-2747 (printed)\hfill {\lnum\number\count0}\else\ifodd\count0
\def\\{ }\ifx\theshorttitle\relax \thetitle \else\theshorttitle\fi\hfill
{\lnum\number\count0}\else\def\\{ and }{\lnum\number\count0}
\hfill\ifx\theshortauthors\relax 
\theauthors\else\theshortauthors\fi\fi\fi}}\def\@evenhead{\@oddhead}
\def\@oddfoot{\small\lfoot\ifnum\count0=\startpage\copyright\ \gtp\hfill\else
\agt, Volume \thevolumenumber\ (\thevolumeyear)\hfill\fi}
\def\@evenfoot{\@oddfoot}
\makeatother
\fi
\let\maketitlepage\makeagttitle
\let\makeshorttitle\maketitlepage
\let\maketitle\maketitlepage

\newwrite\gtoutfile
\long\gdef\makeheadfile{  
{\def\\{, }\def\s{ }
\immediate\openout\gtoutfile head.xxx
\immediate\write\gtoutfile{To: math@arxiv.org}
\immediate\write\gtoutfile{Subject: put OR rep NNNNN:ppppp}
\immediate\write\gtoutfile{--text follows this line--}
\immediate\write\gtoutfile{Proxy-for: \ifx\theasciiauthors\relax
\theauthors\else\theasciiauthors\fi\s<\ifx\theasciiemail\relax\theemail\else\theasciiemail\fi>}
\immediate\write\gtoutfile{\noexpand\\}
\immediate\write\gtoutfile{Authors: \ifx\theasciiauthors\relax
\theauthors\else\theasciiauthors\fi}
{\def\\{ }\immediate\write\gtoutfile{Title: \ifx\theasciititle\relax
\thetitle\else\theasciititle\fi}}
\immediate\write\gtoutfile{Subj-class: GT or SG, GR etc}
\immediate\write\gtoutfile{MSC-class: \theprimaryclass\ifx\thesecondaryclass\relax\else, \thesecondaryclass\fi}
\immediate\write\gtoutfile{Journal-ref: Algebr. Geom. Topol. \thevolumenumber\s
(\thevolumeyear) \startpage-\finishpage}
\immediate\write\gtoutfile{Comments: Published by Algebraic and
Geometric Topology at}
\immediate\write\gtoutfile{\s\s\s  http://www.maths.warwick.ac.uk/agt/AGTVol\thevolumenumber/agt-\thevolumenumber-\thepapernumber.abs.html}
\immediate\write\gtoutfile{\noexpand\\}
\immediate\write\gtoutfile{}
\ifx\theasciiabstract\relax
\immediate\write\gtoutfile{\theabstract}\else
\immediate\write\gtoutfile{\theasciiabstract}\fi
\immediate\write\gtoutfile{}
\immediate\write\gtoutfile{\noexpand\\}
\immediate\write\gtoutfile{}
\immediate\closeout\gtoutfile}}  

\def\maketitlepage{\makeagttitle\makeheadfile}
\let\makeshorttitle\maketitlepage
\let\maketitle\maketitlepage

\def\ifplaintex{\expandafter\ifx\csname documentclass\endcsname\relax}

\def\gtp{{\mathsurround=0pt\it $\cal G\mskip-2mu$eometry \&\ 
$\cal T\!\!$opology $\cal P\!$ublications}}  

\def\recd{{\small Received:\qua\receiveddate\ifx\reviseddate\relax
\else\qquad Revised:\qua\reviseddate\fi\par}} 

\def\lognumber#1{\def\thelognumber{#1}}
\def\volumenumber#1{\def\thevolumenumber{#1}}
\def\volumeyear#1{\def\thevolumeyear{#1}}
\def\papernumber#1{\def\thepapernumber{#1}}
\def\pagenumbers#1#2{\def\startpage{#1}\def\finishpage{#2}}
\def\published#1{\def\publishdate{#1}}

\def\received#1{\def\receiveddate{#1}}
\def\revised#1{\def\reviseddate{#1}}
\def\accepted#1{\def\accepteddate{#1}}

\long\def\asciiabstract#1{\long\def\theasciiabstract{#1}}
\def\asciikeywords#1{\def\theasciikeywords{#1}}

\let\\\par\let\thelognumber\relax\let\thevolumenumber\relax
\let\thepapernumber\relax\let\thevolumeyear\relax\let\startpage\relax
\let\finishpage\relax\let\publishdate\relax\let\receiveddate\relax
\let\reviseddate\relax\let\accepteddate\relax\let\theasciititle\relax
\let\theasciiauthors\relax
\let\theasciiabstract\relax\let\theasciikeywords\relax

\let\theasciiemail\relax

\ifplaintex
\font\logobig=cmssbx10 scaled 3836
\font\logomed=cmssbx10 scaled 2557
\else
\font\logobig=cmssbx10 scaled 4200
\font\logomed=cmssbx10 scaled 2800
\fi

\long\def\makeagttitle{   
\count0=\startpage
\agt\hfill      
\hbox to 45truept{\vbox to 0pt{\vglue -13truept{\logomed A\kern -.37em{\logobig 
T}\kern -.38em G}\vss}\hss}
\break
{\small Volume \thevolumenumber\ (\thevolumeyear)
\startpage--\finishpage\nl
Published: \publishdate}

\vglue .25truein

{\parskip=0pt\leftskip 0pt plus
1fil\def\\{\par\smallskip}{\Large\bf\thetitle}\par\medskip} \vglue
0.05truein

{\parskip=0pt\leftskip 0pt plus 1fil\def\\{\par}{\sc\theauthors}
\par\medskip}%
 
\vglue 0.03truein 

{\small\leftskip 25truept\rightskip 25truept{\bf Abstract}\stdspace\theabstract

{\bf AMS Classification}\stdspace\theprimaryclass
\ifx\thesecondaryclass\relax\else; \thesecondaryclass\fi\par
{\bf Keywords}\stdspace \thekeywords\par}\vglue 7truept

}   

\ifplaintex
\font\phead=cmsl9 scaled 950
\font\pnum=cmbx10 scaled 913
\font\pfoot=cmsl9 scaled 950
\headline{\vbox to 0pt{\vskip -4.5mm\line{\small\phead\ifnum
\count0=\startpage ISSN 1472-2739 (on-line) 1472-2747 (printed)
\hfill {\pnum\folio}\else\ifodd\count0\def\\{ }%
\ifx\theshorttitle\relax\thetitle\else\theshorttitle\fi\hfill{\pnum\folio}
\else\def\\{ and }{\pnum\folio}\hfill\ifx\theshortauthors\relax\theauthors
\else\theshortauthors\fi\fi\fi}\vss}}
\footline{\vbox to 0pt{\vglue 0mm\line{\small\pfoot\ifnum\count0=\startpage
\copyright\ \gtp\hfill\else
\agt, Volume \thevolumenumber\ (\thevolumeyear)\hfill\fi}\vss}}
\else
\font=cmsl9 scaled 1050
\font\lnum=cmbx10 
\font=cmsl9 scaled 1050
\makeatletter
\def\@oddhead{{\small\lhead\ifnum\count0=\startpage ISSN 1472-2739 
(on-line) 1472-2747 (printed)\hfill {\lnum\number\count0}\else\ifodd\count0
\def\\{ }\ifx\theshorttitle\relax \thetitle \else\theshorttitle\fi\hfill
{\lnum\number\count0}\else\def\\{ and }{\lnum\number\count0}
\hfill\ifx\theshortauthors\relax 
\theauthors\else\theshortauthors\fi\fi\fi}}\def\@evenhead{\@oddhead}
\def\@oddfoot{\small\lfoot\ifnum\count0=\startpage\copyright\ \gtp\hfill\else
\agt, Volume \thevolumenumber\ (\thevolumeyear)\hfill\fi}
\def\@evenfoot{\@oddfoot}
\makeatother
\fi
\let\maketitlepage\makeagttitle
\let\makeshorttitle\maketitlepage
\let\maketitle\maketitlepage

\newwrite\gtoutfile
\long\gdef\makeheadfile{  
{\def\\{, }\def\s{ }
\immediate\openout\gtoutfile head.xxx
\immediate\write\gtoutfile{To: math@arxiv.org}
\immediate\write\gtoutfile{Subject: put OR rep NNNNN:ppppp}
\immediate\write\gtoutfile{--text follows this line--}
\immediate\write\gtoutfile{Proxy-for: \ifx\theasciiauthors\relax
\theauthors\else\theasciiauthors\fi\s<\ifx\theasciiemail\relax\theemail\else\theasciiemail\fi>}
\immediate\write\gtoutfile{\noexpand\\}
\immediate\write\gtoutfile{Authors: \ifx\theasciiauthors\relax
\theauthors\else\theasciiauthors\fi}
{\def\\{ }\immediate\write\gtoutfile{Title: \ifx\theasciititle\relax
\thetitle\else\theasciititle\fi}}
\immediate\write\gtoutfile{Subj-class: GT or SG, GR etc}
\immediate\write\gtoutfile{MSC-class: \theprimaryclass\ifx\thesecondaryclass\relax\else, \thesecondaryclass\fi}
\immediate\write\gtoutfile{Journal-ref: Algebr. Geom. Topol. \thevolumenumber\s
(\thevolumeyear) \startpage-\finishpage}
\immediate\write\gtoutfile{Comments: Published by Algebraic and
Geometric Topology at}
\immediate\write\gtoutfile{\s\s\s  http://www.maths.warwick.ac.uk/agt/AGTVol\thevolumenumber/agt-\thevolumenumber-\thepapernumber.abs.html}
\immediate\write\gtoutfile{\noexpand\\}
\immediate\write\gtoutfile{}
\ifx\theasciiabstract\relax
\immediate\write\gtoutfile{\theabstract}\else
\immediate\write\gtoutfile{\theasciiabstract}\fi
\immediate\write\gtoutfile{}
\immediate\write\gtoutfile{\noexpand\\}
\immediate\write\gtoutfile{}
\immediate\closeout\gtoutfile}}  

\def\maketitlepage{\makeagttitle\makeheadfile}
\let\makeshorttitle\maketitlepage
\let\maketitle\maketitlepage

\lognumber{38}
\volumenumber{1}
\volumeyear{2001}
\papernumber{38}
\pagenumbers{763}{790}
\received{18 April 2001}
\revised{27 July 2001}
\accepted{3 October 2001}
\published{11 December 2001}

\usepackage{amsmath}
\usepackage{amsfonts}
\usepackage {epsf}

\newtheoremstyle{plain}{14pt plus6.3pt minus6.3pt}{7.4pt plus3pt minus3pt}%
{\sl}{}{\bf}{}{0.75em}{\thmnumber{#2 }\thmname{#1}\thmnote{\rm\stdspace(#3)}}
\newtheoremstyle{definition}{14pt plus6.3pt minus6.3pt}{7.4pt plus3pt minus3pt}%
{\rm}{}{\bf}{}{0.75em}{\thmnumber{#2 }\thmname{#1}\thmnote{\sl\stdspace#3}}
\newtheoremstyle{remark}{14pt plus6.3pt minus6.3pt}{7.4pt plus3pt minus3pt}%
{\rm}{}{\bf}{}{0.75em}{\thmnumber{#2 }\thmname{#1}\thmnote{\sl\stdspace#3}}

\theoremstyle{plain}
\newtheorem{thm}{Theorem}[section]

\newtheorem{lem}[thm]{Lemma}
\newtheorem{cor}[thm]{Corollary}
\newtheorem{conj}[thm]{Conjecture}

\newtheorem{prop}[thm]{Proposition}

\theoremstyle{definition}

\newtheorem{ex}[thm]{Example}

\theoremstyle{remark}
\newtheorem*{rem}{Remark}

\newtheorem*{claim}{Claim}

\newtheorem*{rems}{Remarks}

\newcommand{\G}{{\Gamma}}
\newcommand{\D}{{\Delta}}
\renewcommand{\O}{{\Omega}}
\newcommand{\bdy}{\partial}
\newcommand{\del}{\partial}
\newcommand{\F}{\mathcal F}

\title{Genus two 3--manifolds are built from\\handle number one pieces}
\author{Eric Sedgwick}
\address{DePaul University, Department of Computer Science\\243 S Wabash Ave, Chicago, IL 60604, USA }
\email{esedgwick@cs.depaul.edu}

\usepackage{amssymb}

\begin{document}

\begin{abstract}
Let $M$ be a closed, irreducible, genus two 3--manifold, and $\F$
a maximal collection of pairwise disjoint, closed, orientable,
incompressible surfaces embedded in $M$.  Then each component
manifold $M_i$ of $M-\F$ has handle number at most one, i.e.\
admits a Heegaard splitting obtained by attaching a single
1--handle to one or two components of $\bdy M_i$.  This result
also holds for a decomposition of $M$ along a maximal collection
of incompressible tori.
\end{abstract}
\asciiabstract{
Let M be a closed, irreducible, genus two 3-manifold, and F
a maximal collection of pairwise disjoint, closed, orientable,
incompressible surfaces embedded in M.  Then each component
manifold M_i of M-F has handle number at most one, i.e.
admits a Heegaard splitting obtained by attaching a single
1-handle to one or two components of boundary M_i.  This result
also holds for a decomposition of M along a maximal collection
of incompressible tori.}


\asciikeywords {3-manifold, Heegaard splitting, incompressible
surface}

\maketitle

\section{Introduction}

Throughout this paper, all surfaces and 3--manifolds will be taken
to be compact and orientable.   Suppose a 3-manifold $M$ contains
an essential 2--sphere.  The Haken lemma \cite{hakenlemma} tells
us that each Heegaard surface for $M$ intersects some essential
2--sphere in a single essential circle (see also \cite{jacobook}).
As a consequence of this and the uniqueness of prime
decompositions of 3-manifolds, Heegaard genus is additive under
connected sum, $$g(M_1 \# \cdots \# M_n) = g(M_1) + \cdots +
g(M_n),$$ where $g(M)$ denotes the Heegaard genus of the manifold
$M$.

How does Heegaard genus behave under decompositions of an
irreducible manifold along incompressible surfaces? Clearly, we do
not expect additivity of genus as before. Suppose that $M$
contains an embedded, incompressible surface $F$ that separates
$M$ into two components $M_1$ and $M_2$. The genus of the two
component manifolds must be greater than the genus of their
boundary component, $g(M_i) > g(F), i=1,2$.  This is particularly
relevant in light of the examples of Eudave-Mu\~noz
\cite{munoz-t1}, tunnel number one knots whose exteriors contain
incompressible surfaces of arbitrarily high genus.  (An
appropriate Dehn surgery on such a knot results in a closed genus
two manifold with an arbitrarily high genus incompressible
surface).

However,  we can build a Heegaard splitting for $M$ from Heegaard
splittings of the components $M_1$ and $M_2$.  If done in an
efficient manner, see for example \cite{s-t-2}, this yields an
upper bound on the genus of $M$, $$g(M) \leq g(M_1) + g(M_2) -
g(F).$$
Upper bounds on the genus of the component manifolds (lower bounds
on $g(M)$) are  more difficult, and not even possible without
additional assumptions.   Consider the examples of Kobayashi
\cite{kob}, knots whose tunnel numbers degenerate arbitrarily
under connected sum (decomposition along an annulus). Again, an
appropriate Dehn surgery will yield a closed manifold containing an
incompressible torus, and after cutting along the torus, the
component manifolds have genus arbitrarily higher than that of the
closed manifold.  In contrast, Schultens \cite{schultens-tunnel}
has demonstrated that for tunnel numbers, this phenomenon cannot occur in the absence
of additional incompressible surfaces. We are led to adding
the assumption that the closed manifold should be cut along
a maximal  embedded collection of incompressible surfaces (a slightly weaker assumption will suffice, see the definition of
a complete collection of surfaces in the next section).

While it is true that the spine of a Heegaard splitting for $M$
induces Heegaard splittings of the component manifolds, see Figure
\ref{f-setup} and Section \ref{minimal}, the intersection between
the Heegaard spine and incompressible surfaces could potentially
be very complicated, and almost certainly depends on the genus of
the incompressible surfaces. One approach to constructing upper
bounds of the genus of the component manifolds is to bound the
complexity of this intersection in terms of the genera of the
incompressible surfaces and the Heegaard spine. This is the
approach used by Johannson in \cite{johannson-book}.

In this paper, we adopt a different approach.  Using ideas of
Scharlemann and Thompson \cite{s-t}, we arrange the spine of the
Heegaard splitting to intersect the collection of surfaces
minimally. It is not hard to see that the induced Heegaard
splitting of the component manifolds is weakly reducible. We then
prove a generalization of a result of Casson and Gordon
\cite{casson-gordon} to manifolds with boundary. (A similar
theorem was proven by Lustig and Moriah \cite{lm}.) A somewhat
simplified version follows:

\begin{thm}
Let $M$ be an irreducible 3--manifold and $M=C_1 \cup_H C_2$ a
weakly reducible Heegaard splitting of $M$.  Then either $M$
contains a closed, non-peripheral incompressible surface, or the
splitting is not of minimal genus.
\end{thm}

This theorem allows us to make use of the assumption that we have
taken the collection of surfaces to be complete. Additional
refinements to this result show that for many of the component
manifolds, the induced Heegaard splitting can be compressed to one
that is induced by a single arc attached to the boundary.

\begin{thm}
\label{t-main} Let $M$ be a closed, irreducible 3--manifold and
$\F$ a complete collection of surfaces for $M$. If
$\overline{M-N(\F)}$ has $n$ component manifolds, then at least $n
+ 2 - g(M)$ of these components have handle number at most 1.
\end{thm}

We give definitions in the next section, but, note here that a
complete collection of surfaces applies both to maximal
collections of incompressible surfaces and maximal collections of
incompressible tori.  Handle number one means that the component
manifold has a Heegaard splitting that is induced by drilling out
a single arc, this is a generalization of tunnel number one and
the concepts are identical when the manifold has a single boundary
component.   While it is possible that a component manifold has
handle number 0, this will occur only when $M$ fibers over the
circle, or unnecessary parallel copies of some surface occur in
the collection. Handle number 0 implies that the component is a
compression body, in fact a product, since its boundary is
incompressible.  This component is either bounded by disjoint
parallel copies of a surface, or there is a single surface cutting
$M$ into a product, i.e., $M$ fibers over the circle.

Unfortunately, we are unable to draw conclusions about every
component manifold unless the genus of $\G$, hence $g(M)$, is 2.

\begin{cor}
If $M$ is a closed, irreducible genus two 3--manifold and $\F$ is
a complete collection of surfaces, then every component manifold
of $\overline{M-N(\F)}$ has handle number at most 1.
\end{cor}

Although the component manifolds have high genus Heegaard
splittings, the fact that they are handle number one means that
this is due precisely to the fact that these manifolds have
boundary with high genus, and are otherwise very simple in terms
of Heegaard structure. The genus of a handle number one component
manifold is bounded above by $$g(M_i) \leq g(\bdy M_i) + 1,$$
where $g(\bdy M_i)$ is the sum of the genera of the components of
$\bdy M_i$. In some cases, this allows one to precisely compute
the genus of the component manifolds. For example, when $g(M)=2$
and the complete collection $\F$ consists of a single separating
surface $F$, we obtain the equality $g(M_i) = g(F) + 1, i=1,2$. By
contrast, in this case Johannson obtains a bound of $g(M_1) +
g(M_2) \leq 2 g(F) + 10$. Or, if $\F$ is a maximal embedded
collection of tori, the single handle is either attached to one
boundary component or is an arc joining two distinct torus
boundary components. As a corollary, we obtain:

\begin{cor}
If $M$ is genus 2, and $\F$ is a maximal collection of tori then
every component of $\overline{M-N(\F)}$ has genus 2.
\end{cor}

Kobayashi  \cite{kob-genus2} has proven a much stronger result
regarding torus decompositions of genus 2 manifolds.

Most of the techniques presented here apply without assumption on
the genus of $M$.   The exception is Proposition
\ref{single-handle} whose hypothesis on handle number will not be
consistently met when the genus of $M$ is greater than two. If the
hypothesis on handle number can be removed, then a general upper
bound on the handle numbers of the component manifolds is
obtained. (It is likely that one must adopt the assumption that
the collection $\F$ is in fact maximal).  This would yield:

\begin{conj}
Let $M$ be a closed, irreducible 3--manifold and $\F$ a maximal
embedded collection of orientable, incompressible surfaces.  If
$M-N(\F)$ has $n$ components then $$\sum_1^n h(M_i) \leq g(M) + n
- 2,$$ where $h(M_i)$ denotes the handle number of the component
manifold $M_i$.
\end{conj}

\section{Preliminaries}

We give brief definitions of concepts related to Heegaard
splittings, the reader is referred to \cite{scharlemann} for a
more thorough treatment.  Let $S$ be a closed surface, $I=[-1,1]$.
A {\it compression body} $C$ is a 3-manifold obtained by attaching
2--handles  and 3--handles to $S \times I$, where no attachment is
performed along $S \times\{1\}$. The boundary of a compression
body is then viewed as having two parts, $\del_+C$ and $\del_-C$,
where $\del_+C= S \times \{1\}$ and $\del_-C=\del C - \del_+C$.
Alternatively, we may construct a compression body $C$ by
attaching 1--handles to $S \times I$ where all attachments are
performed along $S \times \{1\}$. In this case, $\del_-C = S
\times \{-1\}$, $\del_+C = \del C - \del_-C$, and the 1--handles
are dual to the 2--handles of the former construction. In either
construction, we adopt the convention that every 2--sphere
boundary component of $\del_-C$ is capped off with a ball.  If
$\del_-C = \emptyset$ then $C$ is called a {\it handlebody}. Note
that a handlebody can also be defined as a connected manifold with
boundary that possesses a {\it complete collection of compressing
disks}, a properly embedded collection of disks (the cores of the
2--handles) which cut the handlebody into a disjoint union of
balls.

A {\it Heegaard splitting} is a decomposition of a (closed or
bounded) 3--manifold,  $M'=C_1 \cup_S C_2$, where $C_1$ and $C_2$
are compression bodies with their positive boundaries identified,
$S = \del_+C_1 = \del_+C_2$.  In this case $S$ will be a closed
surface embedded in $M'$ and will be called a {\it Heegaard
surface} for $M'$.  The genus of $M'$ is $$g(M') = \min \{ g(S) |
S \text{ is a Heegaard surface for } M' \}.$$
A Heegaard splitting will be called {\it weakly reducible} if
there are non-empty properly embedded collections of compressing
disks $\D_1 \subset C_1$ and $\D_2 \subset C_2$ so that $\bdy D_1
\cap \bdy \D_2 = \emptyset$ in the Heegaard surface $S$.  If it
exists, the collection $\D_0 = \D_1 \cup \D_2$ is called a {\it
weak reducing system} for the Heegaard splitting.

\begin{figure}[ht!]
{\epsfxsize = 4.5 in \centerline{\epsfbox{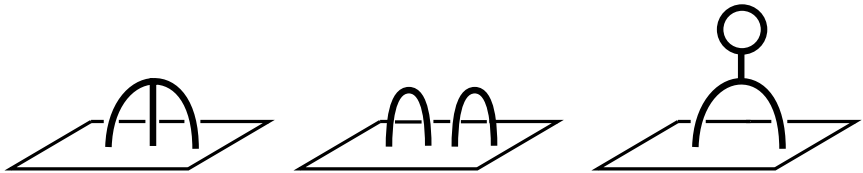}} }
\caption{Graphs with handle number 2.} \label{f-handle}
\end{figure}

If $\G$ is a graph then we will refer to the vertices of valence
1, as the boundary of $\G$, $\bdy \G$.   A graph $\G \subset M'$
will be said to be {\it properly embedded} if it is embedded in
$M$ and $\G \cap \bdy M' = \bdy \G$.   For a properly embedded
graph $\G \subset M'$, we will define the genus of $\G$ to be
$$g(\G) = rank ~ H_1(\G),$$ and define the {\it handle number} of
$\G$ to be $$h(\G) = rank ~ H_1(\G,\bdy \G).$$  Equivalently, the
handle number is the number of edges that need to be removed from
$\G$ so that the resulting graph is empty or a collection of trees
each attached to a boundary component of $M'$ by a single vertex;
or,  $h(\G) = -\chi(\G) + | \bdy \G| =  g(\G) + |\bdy \G| - |\G|$.
Some handle number two graphs are pictured in Figure
\ref{f-handle}.

Typically we will keep track of a Heegaard splitting via a
properly embedded graph in the manifold.  A Heegaard splitting of
a closed manifold $M$ will necessarily consist of two
handlebodies, and in this case, each of the handlebodies is
isotopic to a regular neighborhood of a (non-unique) graph
embedded in the handlebody, hence the manifold. Any such graph
$\G$, for either handlebody, will be called a {\it spine} of the
Heegaard splitting. For bounded manifolds,  Heegaard splittings
come in two varieties, depending on whether or not one of the
compression bodies is actually a handlebody. Correspondingly,
there are two ways that a properly embedded graph can represent a
Heegaard splitting of a bounded manifold.  A {\it tunnel system}
for $M'$ is a properly embedded graph $\G$ so that $\overline{M' -
N(\G)}$ is a handlebody. The {\it tunnel number} of $M'$ is
$$t(M') = \min \{ h(\G) | \G \text{ is a tunnel system for } M'
\}.$$   A {\it handle system} for $M'$ is a properly embedded
graph $\G$ so that $\overline{M'-N(\G)}$ is a compression body $C$
and $\del_-C \subset \del M'$. The {\it handle number} of $M'$ is
$$h(M') = \min \{ h(\G) | \G \text{ is a handle system for } M'
\}.$$  In either case, if $\bdy_1M'$ denotes the boundary
components of $M'$ to which $\G$ is attached, then $\bdy N(\G \cup
\bdy_1M')$ is a Heegaard surface for $M'$.

\begin{figure}[ht!]
{\epsfxsize = 3 in \centerline{\epsfbox{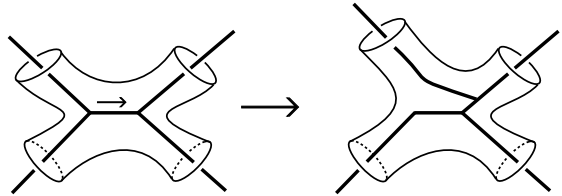}} }
\caption{Edge slides do not change the Heegaard splitting.}
\label{f-eslide}
\end{figure}

Whenever a Heegaard splitting is represented by an embedded graph,
whether a spine, tunnel system, or handle system, then we may
perform slides of edges of the graph along other edges of the
graph without changing the isotopy class of the Heegaard surface,
see Figure \ref{f-eslide}.   Such moves are called {\it edge
slides} or {\it handle slides}.   When working with tunnel or
handle system, we may also slide handles along the boundary of the
manifold without changing the Heegaard splitting.

In the case of a tunnel or handle system, $\G$ will be
slide-equivalent to a collection of $h(\G)$ properly embedded arcs
in $M'$. So $t(M')$ and $h(M')$ should be thought of as the
minimal number of arcs that need to be drilled out of $M'$ so that
the resulting manifold is a handlebody or compression body,
respectively. The handle number is a strict generalization of the
tunnel number and we have $h(M') \leq t(M')$. In general these
quantities are different. For example the exterior of the Hopf
link in $S^3$ is tunnel number one but handle number 0.

A bounded manifold $M'$ will be said to be {\it indecomposable} if
it contains no closed, orientable, non-peripheral incompressible
surface whose genus is either less than or equal to the genus of a
single boundary component of $M'$ or strictly less than the sum of
the genera of two distinct boundary components of $M'$.  Let $\F =
F_1 \cup F_2 \cup \dots,F_k \subset M$ be an embedded collection
of closed orientable incompressible surfaces. A  {\it component
manifold} is a component of the manifold $\overline{M-N(\F)}$. If
$\F$ is an embedded collection of closed, orientable,
incompressible surfaces and each of the component manifolds is
indecomposable, then we say that $\F$ is a {\it complete
collection of surfaces}.

Clearly a maximal embedded collection of orientable,
incompressible surfaces is complete. However, this is not required
for the collection be complete. For example, a maximal embedded
collection of tori in an irreducible manifold is complete as each
of the component manifolds is indecomposable (any additional
surface would have to be genus 1).

\section{Proof of the Main Theorem}
\label{s-main}

In this section we will give an outline of the proof of the main
theorem, Theorem \ref{t-main}.  The proofs of several important
lemmas will be deferred to later sections. Throughout, $M$ will
denote a closed, orientable, irreducible 3--manifold, $\G$ will be
the spine of an irreducible Heegaard splitting of $M$, and $\F$
will be a complete collection of incompressible surfaces.

\begin{figure}[ht!]
{\epsfxsize = 4 in \centerline{\epsfbox{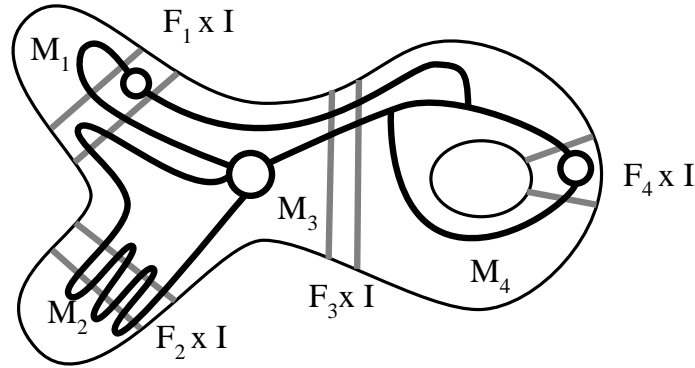}} }
\caption{Letting $\G$ intersect the decomposition minimally.}
\label{f-setup}
\end{figure}

\medskip
{\bf Arranging $\G$ to intersect the decomposition minimally}\qua
Embed in $M$ two parallel copies of each of the incompressible
surfaces in $\F$ and denote this collection by $2~\F$.  If there
are $k$ components of $\F$, $2~\F$ decomposes $M$ into $n+k$
pieces, $k$ {\it product manifolds} $F_j \times I, j=1..k $ and
$n$ {\it component manifolds} denoted $M_i, i=1..n, n < k $,
identical to those obtained by cutting along $\F$.  See Figure
\ref{f-setup}.

Suppose that $\G$ is in general position with respect to $2~\F$
and that we have chosen $\Delta$ to be a complete collection of
compressing disks for the complementary handlebody $\overline{M -
N(\G)}$. The {\it complexity of $(\G,\D)$} is an ordered triple
$(\cdot,\cdot,\cdot)$ of the following quantities:
\begin{enumerate}
\item $\sum h(\G \cap M_i) = $  the sum of handle numbers of the
intersection of the spine $\G$ with each of the component
manifolds $M_i$,
\item $\sum h(\G \cap (F_j \times I)) = $ the sum of the handle
numbers of the intersection of the spine $\G$ with each of the
product manifolds $F_j \times I$,
\item $| \Delta \cap 2~\F |= $ the number of components in the intersection
of the disk collection $\D$ and the surfaces $2~\F$.
\end{enumerate}

Isotoping or manipulating $\G$ by edge--slides  does not change
the isotopy class of the Heegaard surface $\bdy N(\G)$, and we
therefore consider such a spine to be {\it equivalent} to $\G$.
With no loss of generality, we will assume that a spine equivalent
to $\G$ and a complete collection of compressing disks $\D$ have
been chosen to minimize complexity with respect to lexicographic
ordering. Specific properties of the intersection $(\G \cup \D)
\cap 2~\F$ will be developed in Section \ref{minimal}; and are
based on the arguments of Scharlemann and Thompson \cite{s-t}. In
particular we will prove:

\begin{thm}
$\G \cap M'$ is a tunnel system for each product or component
manifold $M'$.
\end{thm}

{\it Proof deferred to Section \ref{minimal}.}

\begin{figure}[ht!]
{\epsfxsize = 3.5 in \centerline{\epsfbox{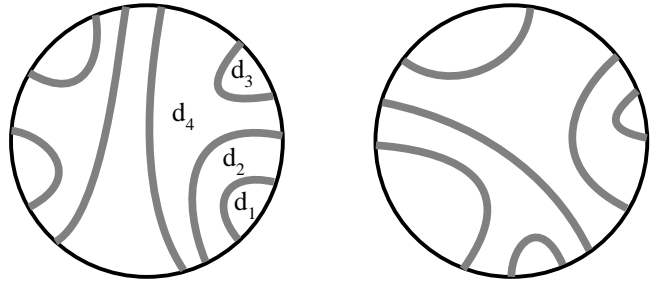}} }
\caption{Labeling subdisks of $\D - 2~\F$.} \label{f-delta}
\end{figure}

\bigskip
{\bf Ordering subdisks of $\D - 2~\F$}\qua By Lemma \ref{t-nocircles}
we know that $\D - 2~\F$ is a collection of disks.  We will
(non-uniquely) label these disks $d_1,\dots,d_m$ according to the
following rules:
\begin{enumerate}
\item Label an outermost disk $d_1$,
\item Assuming that the disks $d_1,\dots,d_{l-1}$ have been labeled,
give the label $d_l$ to a subdisk of $\D - 2~\F$ that is outermost
{\it relative to} the subdisks $d_1,\dots,d_{l-1}$. See Figure
\ref{f-delta}.
\end{enumerate}

Note that each of the subdisks $d_l$ is embedded in some component
or product manifold $M'$.  Moreover, it is a compressing disk for
the handlebody that is the complement of the tunnel system induced
by $\G$, $\overline{M'-N(\G)}$.

Let $\{d_{i_j}\}$ be a subcollection of the disks $\Delta - 2~\F$.
The {\it support of $\{d_{i_j}\}$}, denoted $supp(\{d_{i_j}\})$,
is the sub--graph of $\G$ that is the spine of the handlebody
obtained by maximally compressing $N(\G)$ along compressing disks
which are disjoint from $2~\F$ and disjoint from the boundary of
$\{d_{i_j}\} \subset \bdy N(\G)$ and throwing away any components
which do not meet $\{d_{i_j}\}$.  For each component manifold
$M_i$ let $j$ be the least $j$ so that $d_j \subset M_i$.  The
disk $D_i = d_j$ will be called {\it the relatively outermost
disk} for $M_i$.   The graph $$\O_i = supp(D_i)$$ will be called
{\it the relatively outermost graph} for $M_j$.

\begin{figure}[ht!]
{\epsfxsize = 4.85 in \centerline{\epsfbox{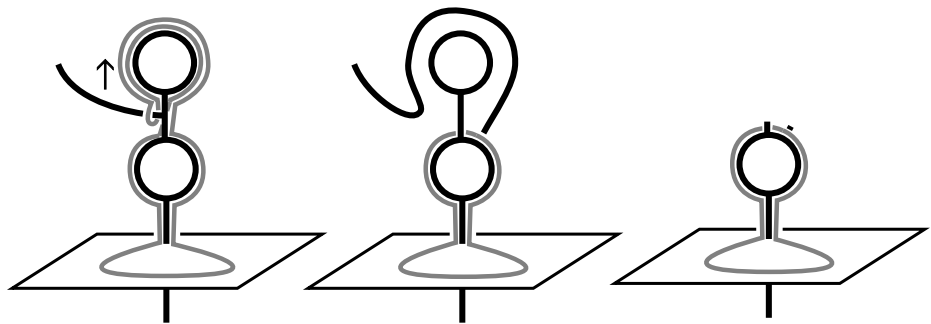}} }
\caption{The support of a disk.} \label{f-supp}
\end{figure}

\vspace{1in}

\begin{rems}
\hspace{5in}
\begin{enumerate}
\item In the definition of support, it may be necessary to
perform handle slides of $\G$ in the interior of some component
manifolds in order to realize the maximal collection of
compressing disks, see Figure \ref{f-supp}.

\item We have chosen a fixed numbering of the subdisks of
$\D - 2~\F$. Thus, the notions of the relatively outermost disk
and the relatively outermost graph for a component manifold are
well defined.

\item We will consider the support of a relatively outermost disk $\O_i=supp(D_i)$
to be a graph that is properly embedded in the component manifold
$M_i$.  We will consider the support of a collection of subdisks
$supp(\{d_{i_j}\})$ to be graph that is embedded in $M$.

\end{enumerate}
\end{rems}

We can reconstruct the spine $\G$ by building a sequence of
graphs, each the support of a larger collection of ordered
subdisks of $\D-2~\F$, $$\Gamma_k = supp(\{d_i\}_{i=1..k}).$$ In
particular, $$\G_m = \G.$$   The relatively outermost graphs for
each component manifold, $\O_1, \dots, \O_n$, will be attached at
some point in building $\G$. Moreover, they are the support of the
relatively outermost disks $D_1,\dots,D_n$, and as we will see
they are attached to the previous graph along all but at most one
of their endpoints.  This gives us a lower bound for the genus of
$\G$ in terms of the handle number of the outermost graphs $\O_i$.

\begin{lem}
Let $\Gamma$ be the spine of an irreducible Heegaard splitting.
Then $$g(\Gamma) \geq \sum_1^n h(\O_i) - n + 2.$$
\end{lem}

{\it Proof deferred until Section \ref{counting}.}
\bigskip

However, it is our aim to develop a lower bound for the genus of
$\G$ in terms of the handle numbers of the component manifolds,
not just the handle numbers of the relatively outermost graphs
$\O_i$.  In a special case ($h(\O_i) = 1$) we will show that
$\O_i$ is in fact a handle system and obtain the desired bound.

\begin{prop}
\label{T-handlesys} If $h(\O_i) = 1$ then $\O_i$ is a handle
system for $M_i$.  In particular,  $h(M_i) \leq 1$.
\end{prop}

{\it Proof deferred to Section \ref{s-weakly-reducible}.}

\begin{rem}
The restriction $h(\O_i) = 1$ in this proposition is what prevents
us from making a more general statement connecting genus to the
sum of handle numbers of the component manifolds.  If $\O_i$ were
always a handle system for the component manifold $M_i$ then we
would obtain the more general inequality $g(\G) \geq \sum h(M_i) -
n + 2$.
\end{rem}

These lemma and proposition prove the main theorem. Let $j \leq n$
be the number of components $M_i$ which have $h(M_i) > 1$. By
Proposition \ref{T-handlesys}, each of the corresponding outer
handle systems $\O_i$ has $h(\O_i) > 1$. By Lemma \ref{T-count} we
have $$g(\Gamma) - 2 \geq \sum h(\O_i) - n \geq n + j - n \geq j
.$$ Therefore, the number of component manifolds with handle
number one is at least $$n - (g(\G)-2) = n+2-g(\G).\eqno{\qed}$$ 

\section{Properties of the Minimal Intersection between $2~\F$ and $\G$}
\label{minimal}

This section is devoted to developing the properties of the
minimal intersection between the Heegaard complex $\G \cup \Delta$
and the incompressible surfaces $2~\F$.  Many of these properties
were either specified in the work of Scharlemann and Thompson
\cite{s-t}, or follow from the same methods. They are included
here, both for the sake of completeness, and due to the fact that
the definition of minimality used here differs from that in
\cite{s-t}. We also apply these properties to characterize the
support of outermost and relatively outermost disks.

Throughout this section, we assume that the spine $\G$ of the
irreducible Heegaard splitting and compressing disks $\D$ for its
complement have been chosen to intersect the surfaces $2~\F$
minimally, as defined in the previous section. However, it is not
necessary to place any restrictions on the surface collection
$\F$.

First we will demonstrate that $\G$ induces Heegaard splittings of
each of the component and product manifolds.

\begin{lem}
\label{t-incomp} Let $F$ be a component of $2~\F$.  Then the
punctured surface $F'=\overline{F-N(\G)}$ is incompressible in the
handlebody $\overline{M - N(\G)}$.
\end{lem}

\begin{figure}[ht!]
{\epsfxsize = 5 in \centerline{\epsfbox{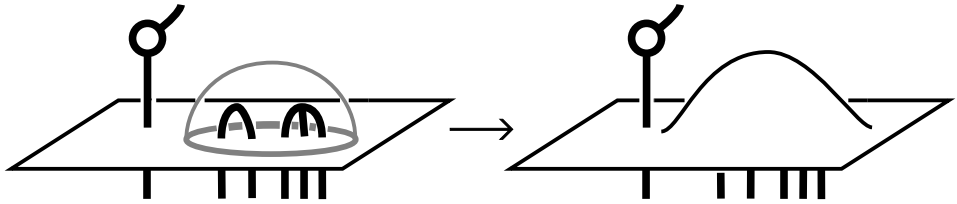}} }
\caption{If $F'$ is compressible the intersection is not minimal.}
\label{f-ball}
\end{figure}

\begin{proof}
If some component of the punctured surface were compressible, then
there would be a compressing disk $D$ for $F'$, a perhaps distinct
component of the punctured surface, embedded in some component or
product manifold $M'$. The boundary of $D$ bounds a disk $D'$ in
$F$. As $M$ is irreducible, $D$ and $D'$ cobound a ball $B$, and
$B$ must be contained in $M'$, for otherwise the incompressible
surface $F$ would lie in the ball $B$.   We can therefore isotope
$F$ through $B$, thereby pushing a portion of $\G \cap M'$ into an
adjacent product or component manifold. See Figure \ref{f-ball}.
Since $\G$ is the spine of an irreducible splitting, by a theorem
of Frohman \cite{frohman}, $B \cap \G$ does not contain any loops
of $\G$, it is merely a collection of trees. As there is no loop
of $\G$ in the ball, this move does not raise the induced handle
number of the adjacent manifold, while it does reduce the handle
number of $\G \cap M'$. This contradicts the minimality of the
intersection.
\end{proof}

\begin{thm}
\label{t-tunnelsys} $\G \cap M'$ is a tunnel system for each
product or component manifold $M'$.
\end{thm}

\begin{proof}
The manifold $\overline{M' - N(\G)}$ is a component of the
handlebody $\overline{M-N(\G)}$ after it is cut along the properly
embedded collection of punctured incompressible surfaces $2 ~\F' =
\overline{2~\F - N(\Gamma)}$. It is well known that when a
handlebody is cut along a collection of incompressible surfaces,
the resulting pieces are handlebodies. So $\overline{M'- N(\G)}$
is a handlebody and $\G \cap M'$ is the corresponding tunnel
system.
\end{proof}

\begin{figure}[ht!]
{\epsfxsize = 4 in \centerline{\epsfbox{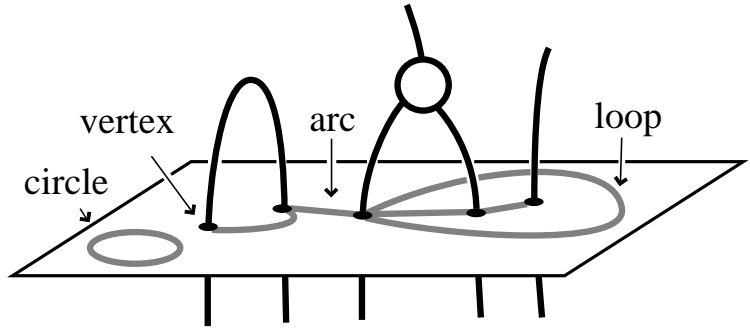}} } \caption{The
intersection of $\G \cup \D$ with $2~\F$ is a graph in $2~\F$.}
\label{f-G}
\end{figure}

The intersection of the 2-complex $\G \cup \D$ with the
incompressible surfaces $2~\F$ is a graph $G \subset 2~\F$. See
Figure \ref{f-G}.   A component of  intersection with the spine,
$\G \cap 2~\F$ is called a {\it vertex}.  Since handlebodies do
not contain closed incompressible surfaces, there is at least one
vertex in each component of $2~\F$. A component of the
intersection with the compressing disks, $\D \cap 2~\F$ is called
a {\it circle} if it is an intersection with the interior of $\D$
and an {\it edge} otherwise.   An edge joining distinct vertices
will be called an {\it arc}  and an edge joining a vertex $v$ to
itself is called a {\it loop based at $v$}.

\begin{lem}
\label{t-nocircles} There are no circles in $G$.
\end{lem}

\begin{proof}
This follows from the minimality of $\D \cap 2~\F$, using an
innermost disk argument and Lemma \ref{t-incomp}.
\end{proof}

\begin{lem}{\rm\cite{s-t}}\qua
There are no isolated vertices (every vertex belongs to some
edge).
\end{lem}

\begin{proof}
If some vertex is isolated then it defines a compressing disk $D$
for the handlebody $N(\G)$ (or the vertex cuts off a tree,
contradicting minimality). Moreover, the boundary of $D$ is
disjoint from the complete collection of disks $\Delta$.  After
compressing the handlebody $\overline{M-N(\G)}$ along $\Delta$ we
obtain a collection of balls, and $\bdy D$ is a loop on the
boundary of one of these balls.  It therefore  also bounds a disk
in the handlebody $\overline{M-N(\G)}$.  This implies that $\G$ is
the spine of a reducible Heegaard splitting.
\end{proof}

We rely heavily on the notion of outermost edges \cite{s-t}. Every
edge $e$ of $G$ separates some disk $D \subset \D$ into two
subdisks, $D_1$ and $D_2$. If one of the subdisks  does not
contain any other edges of $G$ then $e$ is called an {\it
outermost edge} of $G$.  Suppose that an edge $e$ is joined to the
vertex $v$ and that one of the two subdisks $D_1$ or $D_2$ does
not contain an edge of $G$ which is joined to $v$. Then, $e$ is an
{\it outermost edge with respect to} $v$.

Note that by passing to subdisks, every vertex $v$ has some edge
$e$ which is outermost with respect to it.  Also, an outermost
edge is outermost with respect to its vertices (or vertex), but
not (in general) vice-versa.

\begin{lem}{\rm\cite{s-t}}\qua
\label{t-outedge} Let $e$ be an outermost edge with respect to one
of its vertices $v$.   Then $e$ is a loop based at $v$ that is
essential in $2~\F$.
\end{lem}

\begin{figure}[ht!]
{\epsfxsize = 3.5 in \centerline{\epsfbox{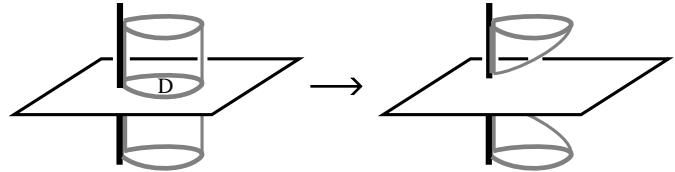}} }
\caption{Boundary compressing  $\Delta$ reduces the number of
intersections.} \label{f-cdisk}
\end{figure}

\begin{proof}
Suppose that $e$ is an arc and joins $v$ to a distinct vertex $w$.
See Figure \ref{f-slide}.  (The edge $e$ may or may not be
outermost for $w$.) The edge $e$ cuts off a disk $D' \subset \D$
which does not contain any edge joined to $v$.  Let $M'$ be the
adjacent manifold into which $D'$ starts, and let $M''$ be
adjacent manifold.  Let $\gamma \subset \G$ denote the handle
containing $v$.

We will now perform a {\it broken edge slide} \cite{s-t} which
shows that the intersection is not minimal.  See Figure
\ref{f-slide}. Add a new vertex to $\gamma$ that lies slightly
into $M'$, this breaks $\gamma$ into two handles, $\gamma_1$ and
$\gamma_2$.  Use the disk $D'$ to guide an edge-slide of
$\gamma_1$, which pulls it back into $M''$. This edge-slide is
permissible precisely because $e$ is outermost for $v$, we did not
ask the handle $\gamma_1$ to slide along itself.    It does not
increase the handle number of $\G \cap M''$, while it strictly
decreases the handle number of $\G \cap M'$ (and possibly others,
if $\gamma_2$ runs through other manifolds). This contradicts the
minimality of the intersection between $\G$ and $2~\F$.

\begin{figure}[ht!]
{\epsfxsize = 4 in \centerline{\epsfbox{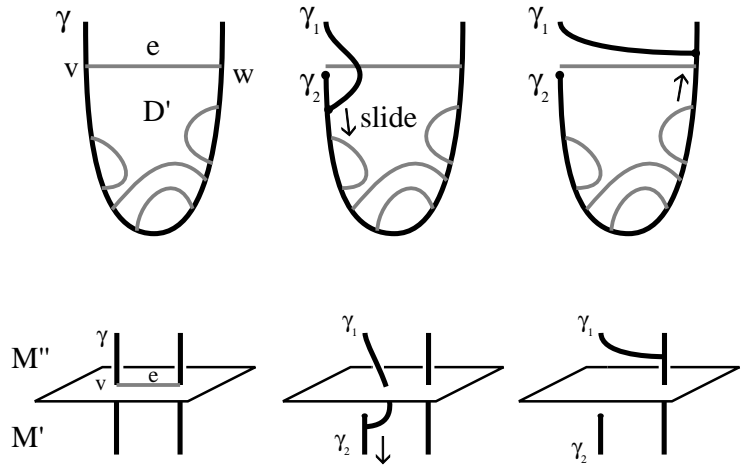}} } \caption{A
broken edge slide. } \label{f-slide}
\end{figure}

We have established that an outermost edge for a vertex must be a
loop.  If it were inessential then we can find an innermost
inessential loop bounding a disk $D$.  If $D$ contains a vertex
$v$ then an outermost edge for $v$ is an arc, contradicting the
previous conclusion of this lemma. If $D$ does not contain a
vertex, then we can reduce the number of intersections of $\D \cap
2~\F$ by boundary compressing $\D$ along $D$. See Figure
\ref{f-cdisk}.
\end{proof}

\begin{lem}
\label{t-outsupp} The support of an outermost disk, $supp(d_j)$ is
connected, has a single boundary vertex, and
$h(supp(d_j))=g(supp(d_j)) > 0$.
\end{lem}

\begin{proof}
There is a single edge $e \subset G$ cutting off the outermost
disk $d_j$ from $\D$. By Lemma \ref{t-outedge}, this edge is an
essential loop in some component $F$ of $2~\F$.  This implies that
$supp(d_j)$ has a single boundary vertex and is connected.  Now,
if $h(supp(d_j)=0$, then the subarc $\alpha = \bdy d_j - e \subset
\bdy N(\G)$ of the boundary of $d_j$ does not cross any
compressing disk of $\G \cap M_i$ other than the disk
corresponding to the vertex. This means that we can perform edge
slides of $\G$ that allow us to pull the arc $\alpha$ back to $F$,
creating an essential circle of intersection in the process. This
is a contradiction, a subdisk of $d_j$ becomes a compressing disk
for $F$, see Figure \ref{f-outer3}.

\begin{figure}[ht!]
{\epsfxsize = 4 in \centerline{\epsfbox{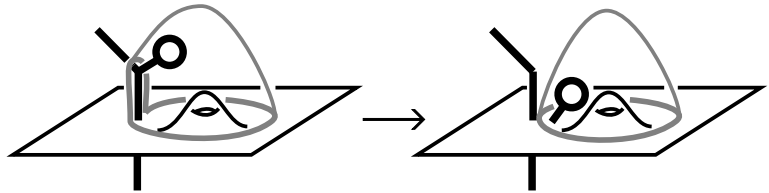}} }
\caption{An outermost disk with handle number 0.} \label{f-outer3}
\end{figure}

Since, $supp(d_j)$ has a single boundary vertex, all of its
handles must be realized by genus, i.e., $g(supp(d_j)) =
h(supp(d_j))$.
\end{proof}

\begin{lem}
\label{t-reloutsupp} The support $\O_i$ of a relatively outermost
disk $D_i$ for a component manifold $M_i$ is connected and has
$h(\O_i) \geq 1$.
\end{lem}

\begin{proof}
We first show that $\O_i$ is connected.  The boundary of the
relatively outermost disk $D_i$ consists of arcs on $\O_i$ and
edges lying in  $2~\F$.  Each arc in $\O_i$ lies in a single
component of $\O_i$. All but at most one of the edges cuts off a
disk which does not return to $M_i$.  Each of these edges is
therefore outermost for its vertices, and by Lemma \ref{t-outedge}
an essential loop in some component $F$ of $2~\F$.  Loops do not
join distinct components of $\O_i$.   This means that $\O_i$ is
connected, for any edge leaving a component there must be an
additional edge that returns to that component, and we have at
most one edge that is not a loop.

\begin{figure}[ht!]
{\epsfxsize = 2 in \centerline{\epsfbox{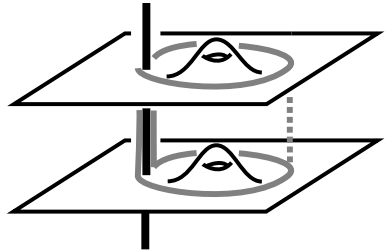}} } \caption{A
relatively outermost disk that joins distinct boundary
components.} \label{f-outer4}
\end{figure}

Now, suppose that $h(\O_i) = 0$.  We know that all but at most one
of the edges is a loop.  While in general it is possible that the
remaining edge $e$ is an edge, this does not occur when $h(\O_i) =
0$.  A single edge implies that the boundary of $D_i$ joins two
distinct vertices in the graph and therefore traverses a handle in
the component manifold.

\begin{figure}[ht!]
{\epsfxsize = 5 in \centerline{\epsfbox{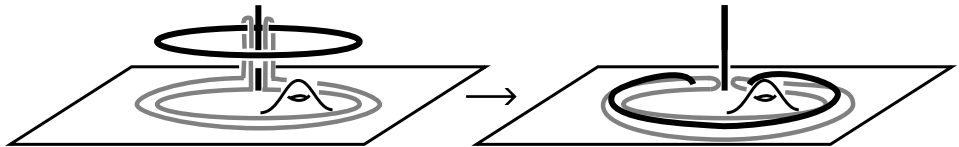}} } \caption{A
relatively outermost disk whose support has handle number 0. }
\label{f-nocap2}
\end{figure}

 Since $h(\O_i)=0$ we may perform edge slides so that a sub disk
of $D_i$ intersects some component $F$ of $2~\F$ in a circle that
bounds a disk in $M_i$.  See Figure \ref{f-nocap2}. This may raise
the handle number of an adjacent product manifold.  Since $F$ is
incompressible, the boundary of this disk bounds a disk in $F$,
the two disks bound a ball, and as in Lemma \ref{t-incomp} we can
perform an isotopy of the graph that reduces the induced handle
number of the component manifold $M_i$. This contradicts
minimality of the intersection of $\G$ and $2~\F$.
\end{proof}

\begin{figure}[ht!]
{\epsfxsize = 2.5 in \centerline{\epsfbox{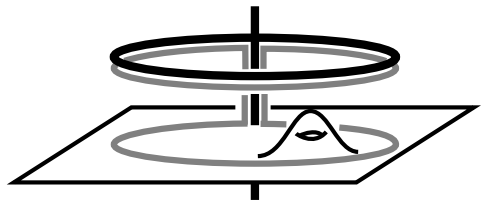}} } \caption{A
handle that is parallel to a component of $\F$.}
\label{f-parallel}
\end{figure}

There is one situation contradicting minimality that cannot be
detected from the intersection of $\G$ and $2~\F$ and the
knowledge that an edge is outermost. It is possible that there is
a loop based at a vertex $v$ that cuts off a disk lying in a
component manifold which runs along a handle exactly once, see for
example Figure \ref{f-parallel}.  In this case, the handle can be
slid into the product manifold reducing complexity. This is also
the motivation for working with the collection $2~\F$ instead of
$\F$ and choosing our definition of complexity.  If we were
working with a single copy, $\F$, this move would not decrease
complexity, it raises the induced handle number of the adjacent
component manifold. This situation will be detected by using the
machinery of Casson and Gordon \cite{casson-gordon} and is
analyzed in Section \ref{s-weakly-reducible}.

\section{Estimating the Genus of $\G$}
\label{counting}

The setup for this section is that of the proof of the main
theorem:  $\G$ is the spine of an irreducible Heegaard splitting,
$\D$ is a complete collection of compressing disk for
$\overline{M-N(\G)}$, both chosen to intersect $2~\F$ minimally;
and $\O_1,\dots,\O_n$ are the support of relatively outermost
disks, $D_1,\dots,D_n$, for the component manifolds,
$M_1,\dots,M_n.$ We demonstrate that the sum of the handle numbers
of the supports gives us a lower bound on the genus of $\G$.

\begin{lem}
\label{T-count} $$g(\Gamma) \geq \sum_1^n h(\O_i) - n + 2$$
\end{lem}

\proof
Recall that we have defined  $$\Gamma_k =
supp(\{d_i\}_{i=1..k}),$$ where $d_1,\dots,d_m$ is an outermost
ordering of the subdisks of $\D - 2~\F$.  The proof is an
inductive one, demonstrating that when $\G_k \subsetneq \G$,
\begin{gather}
\label{E-complexity} g(\Gamma_k) - |\G_k| \geq \sum_{\O_i \subset
\G_k} (h(\O_i)-1),
\end{gather}
where $|\G_k|$ denotes the number of components of $\G_k$.  We
then analyze the final attachment, when $\G_k=\G$.   Note that at
each stage we are attaching some portion of the spine $\G$; the
right hand side of the inequality can only increase when this
portion is actually the support of the relatively outermost disk
$\O_i$ for some component manifold $M_i$.

Let $k=1$. The graph $\G_1$ is the support of the outermost disk
$d_1$ which is embedded in either a component or product manifold
$M'$. By Lemma \ref{t-outsupp}, $\G_1$ is connected, has a single
boundary vertex, and has positive genus.  This means that $g(\G_1)
= h(\G_1)$.  If $M'$ is a component manifold we have $g(\G_1)-1
\geq h(\O_1)-1$ and if $M'$ is a product manifold we have
$g(\G_1)-1 \geq 0$.  This establishes Inequality
\ref{E-complexity} for $k=1$.

Now, suppose that $k > 1$, $\G_k \subsetneq \G$, and that
$\G_{k-1}$ satisfies the inductive hypothesis. If $d_k$ is not an
relatively outermost disk for a component manifold, then we merely
need to observe that the left-hand side of Inequality
\ref{E-complexity} does not decrease when we attach $supp(d_k)$.
It will decrease only if the number of components increases, which
means that some component of $supp(d_k)$ is not attached to
$\G_{k-1}$. But, this happens only if $d_k$ is an outermost disk,
in which case $supp(d_k)$ has a single component and there is an
increase of genus to compensate for the additional component.

We are left in the case that $d_k$ is a relatively outermost disk
$D_i$ for some component manifold $M_i$, then $supp(d_k)=\O_i$. By
Lemma \ref{t-reloutsupp}, $\O_i$ is connected.  Again, if $\O_i$
is not actually attached to $\G_{k-1}$,  an additional component
is added, but then $d_k$ is actually an outermost disk,
$supp(d_k)$ is connected, has positive genus, and
$h(\O_i)-1=g(\O_i)-1$ is added to both sides.

If $d_k$ is a relatively outermost disk $D_i$ for $M_i$, but not
absolutely outermost (for example that in Figure \ref{f-outer4}),
then all but at most one boundary vertex of $\O_i$ is attached to
$\G_{k-1}$. As noted in the proof of Lemma \ref{t-reloutsupp}, all
but at most one of the vertices of $\O_i$, has an outermost loop
in $G$ attached to it that cuts off a subdisk of $\D$ containing
only disks with labels $d_i$, where $i<k$. Each such edge of the
disk is attached to a disk with strictly smaller labels,  so all
but one boundary vertex is attached to $\G_{k-1}$.

\begin{figure}[ht!]
{\epsfxsize = 2.5 in \centerline{\epsfbox{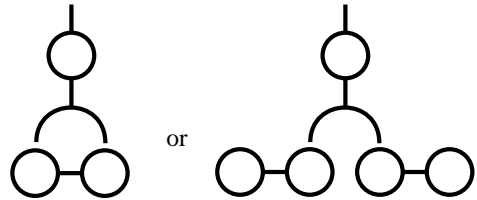}} }
\caption{Attaching  $\O_i$, a relatively outermost graph with
handle number 3, adds 2 to $g(\G_{k-1})-|\G|$.} \label{f-attach}
\end{figure}

So, for all but the first vertex attached, each attached vertex
adds to the genus by one or reduces the number of components by 1,
see Figure \ref{f-attach}. Moreover, any genus of $\O_i$ is added
to the genus of $\G_{k-1}$. We have added $g(\O_i) + |\bdy \O_i| -
2$ to the left-hand side of \ref{E-complexity}, this is the same
as $h(\O_i) - 1$, which is added to the right side.

\begin{figure}[ht!]
{\epsfxsize = 2.5 in \centerline{\epsfbox{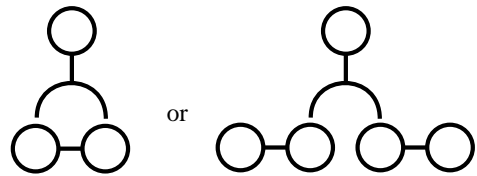}} }
\caption{The final attachment.} \label{f-fattch}
\end{figure}

A similar analysis pertains for the final attachment, when $\G_k =
\G$. However, in this case {\it every} vertex of $supp(d_k)$ is
attached to $\G_{k-1}$, for there can be no unattached vertices,
see Figure \ref{f-fattch}. When $supp(d_k)$ is not some $\O_i$
this adds at least 1 to left side of Inequality \ref{E-complexity}
and nothing to the right hand side. When $supp(d_k) = \O_i$ for
some $i$, this adds $h(\O_i)$ to the left hand side and
$h(\O_i)-1$ to the right side. In either case, the inequality will
still hold even if we add an additional 1 to the right side. This
yields $$g(\G) - |\G| \geq \sum_1^n (h(\O_i) - 1) + 1.$$ Since
$\G$ is connected, we have $$g(\G) \geq \sum_1^n h(\O_i) - n +2
.\eqno{\qed}$$

\section{Weakly Reducible Heegaard splittings of Manifolds with Boundary}
\label{s-weakly-reducible}

In \cite{casson-gordon} Casson and Gordon introduced the notion of
a weakly reducible Heegaard splitting of a closed 3--manifold, and
showed that such a splitting is either reducible or the manifold
contains an incompressible surface.  We first state and prove an
extension of their theorem to manifolds with boundary; a similar
theorem was proven by Lustig and Moriah \cite{lm}.   We will then
apply these techniques to the Heegaard splittings of the component
manifolds $M_i$ that are induced by the Heegaard spine $\G$. These
splittings are typically weakly reducible.

First, we introduce some notation.  Suppose that $H \subset M$ is
a closed embedded surface and $\D \subset M$ is an embedded
collection of disks so that $\D \cap H = \bdy \D$. Let
$\sigma(H,D)$ denote the surface obtained by performing an ambient
1-surgery of $H$ along $D$ (i.e., compression).  We use the notion
of complexity introduced in \cite{casson-gordon}, the {\it
complexity} of a surface is defined to be, $$c(surface) = \sum (1
- \chi(S)),$$ where the sum is taken over all non-sphere
components of the surface $S$.  Note that if $D$ is a single disk
with essential boundary then
\begin{gather*}
c(\sigma(H,D)) =
  \begin{cases}
    c(H) -1 & \text{ if $D$ is separating or compresses a torus} \\
            & \text{ component of $H$, and}, \\
    c(H) - 2 & \text{otherwise}.
  \end{cases}
\end{gather*}

\begin{thm}
\label{weakly-reducible} Let $M$ be a compact, orientable,
irreducible 3-manifold and $M=C_1 \cup_H C_2$ a Heegaard splitting
of $M$.  If $\D_0 = \D_1 \cup \D_2$ is a weak reducing system for
the Heegaard splitting then either
\begin{enumerate}
\item $M$ contains an orientable, non-peripheral incompressible
surface $S$ so that $c(S)$ is less than or equal to the complexity
of $\sigma(H,\D_0)$, or
\item there is an embedded collection of disks $\widehat{\D_1} \subset C_1$
so that $\D_1 \subset \widehat{\D_1}$ and some component of
$\sigma(H,\widehat{\D_1})$ is a Heegaard surface for $M$, or
\item there is an embedded collection of disks $\widehat{\D_2} \subset C_2$
so that $\D_2 \subset \widehat{\D_2}$ and some component of
$\sigma(H,\widehat{\D_2})$ is a Heegaard surface for $M$.
\end{enumerate}
In particular, conclusions (2) and (3) imply that $H$ is not of
minimal genus.
\end{thm}

\begin{proof}
Define the complexity of a weak reducing system  $\D_0 = \D_1 \cup
\D_2$ for the Heegaard splitting $C_1 \cup_H C_2$ to be $$c(\D_0)=
c(\sigma(H,\D_0)).$$  Let the surfaces $H_i = \sigma(H,\D_i),
i=0,1,2$, be obtained by compressing $H$ along the corresponding
disk collections.  See the schematic in Figure \ref{f-weakred}, it
is essential to understanding the arguments of this section.  Note
that the surface $H_0$ separates $M$ into two components, denote
these by $X_1$ and $X_2$.   If we compress the compression body
$C_1$ along the disk system $\D_1$ we obtain a compression body
$Y_1$ which we will think of as sitting slightly inside $X_1$. Its
complement $X_1-Y_1$ can be thought of as $(H_1 \times I) \cup
N(\D_2)$ - a product with 2-handles attached, and is therefore a
compression body. Symmetrically, we also have that $X_2 -Y_2$ is a
compression body. Thus, the surfaces $H_1=\bdy Y_1 $ and $H_2=\bdy
Y_2$ are Heegaard surfaces for the (possibly disconnected)
manifolds $X_1$ and $X_2$, respectively.

\begin{figure}[ht!]
{\epsfxsize = 5.35 in \centerline{\epsfbox{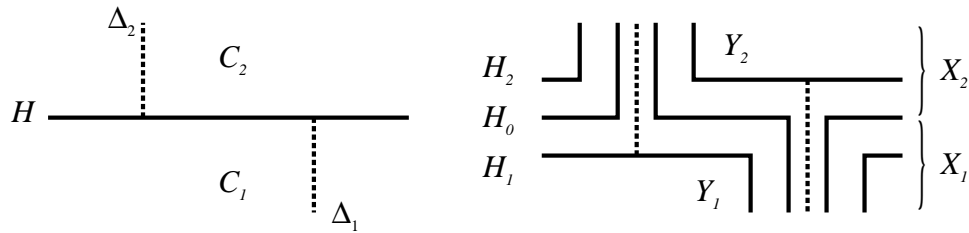}} }
\caption{Compressing a Heegaard surface along a weak reducing
system.} \label{f-weakred}
\end{figure}

Suppose that some  positive genus component of $H_0$ is
compressible, say into $X_1$. The compressing disk $D$ for $H_0$
is a boundary reducing disk for the manifold $X_1$. As $H_1$ is a
Heegaard surface for $X_1$, the Haken lemma (see
\cite{casson-gordon}) implies that we may isotope $D$ to intersect
$H_1$ in a single circle.  It also says that we may choose a new
collection of compressing disks $\D_2'$ for $X_1-Y_1$, hence for
$C_2$, which is disjoint from $D$.  The collection $\D_0' = (\D_1
\cup D) \cup \D_2'$ is a weak reducing system with lower
complexity than $\D_0$ because we have compressed $H$ along an
additional disk.  A symmetric phenomenon occurs if $H_0$ is
compressible into $X_2$.

In fact, we may continue to compress $H_0$, finding new disk
collections of strictly decreasing complexity, until each
component of $H_0$ is a 2-sphere or incompressible surface. Denote
the final weak reducing system by $\D_0'$, and the corresponding
surfaces and sub-manifolds indicated in Figure \ref{f-weakred} by
$H_1',X_1',Y_1',\dots$, etc..  Now $\D_0'$ may or may not contain
the original disk collections $\D_1$ and $\D_2$. However, the
compression bodies $Y_1'$ and $Y_2'$ are obtained by compressing
the compression bodies $Y_1$ and $Y_2$. These in turn were
obtained by compressing $C_1$ and $C_2$ along the original
collections $\D_1$ and $\D_2$. So we may also think of $Y_1'$ and
$Y_2'$ as being obtained by compressing $C_1$ and $C_2$ along a
collection of disks $\widehat{\D_1} \subset C_1$ and
$\widehat{\D_2} \subset C_2$ where $\D_1 \subset \widehat{\D_1}$
and $\D_2 \subset \widehat{\D_2}$. In general $\widehat{\D_1}$ and
$\widehat{\D_2}$ do not have disjoint boundary on $H$ and cannot
be taken to be part of a weak reducing system.

If some component $S$ of $H_0'$ is an incompressible and
non-peripheral surface, then we have conclusion (1) of the
theorem.  Moreover, we have that $c(S) \leq c(\sigma(H,\D_0))$,
for $S$  is a component of $H_0'$ which was obtained by
compressing $H_0 = c(\sigma(H,\D_0))$.

\begin{figure}[ht!]
{\epsfxsize = 2 in \centerline{\epsfbox{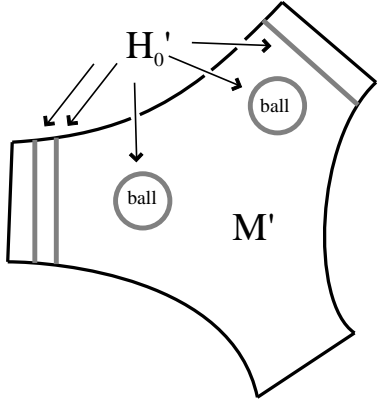}} } \caption{A
component $M'$ is essentially the same as $M$.} \label{f-mprime}
\end{figure}

We therefore assume that each component of $H_0'$ is a 2-sphere
bounding a ball ($M$ is irreducible) or a peripheral
incompressible surface.  See Figure \ref{f-mprime}. Then some
component of $M - N(H_0')$, call it $M'$, is a copy of $M$ perhaps
with some balls and product neighborhoods of  boundary components
of $M$ removed. Since $H_0'$ separates $M$ into $X_1'$ and $X_2'$,
$M'$ must actually be a component of either $X_1'$ or $X_2'$, say
$X_1'$. Recall that the surface $H_1'$ is a Heegaard surface for
$X_1'$. This means that some component  $H'' \subset H_1'$ is a
Heegaard surface for $M'$. In fact, $H''$ is also a Heegaard
surface for $M$; filling in balls and product neighborhoods of the
appropriate boundary components does not change the property that
$H''$ bounds compression bodies to both sides. Moreover, the
Heegaard surface $H''$ is a component of the boundary of $Y_1'$
and it follows from our earlier remarks, that it is a component of
the surface $\sigma(H,\widehat{\D_1})$, where $\D_1 \subset
\widehat{\D_1}$. Symmetrically, if $M' \subset X_2'$ then the
Heegaard surface $H''$ is a component of the surface
$\sigma(H,\widehat{\D_2})$, where $\D_2 \subset \widehat{\D_2}$.
\end{proof}

There is one major difference between the case of closed manifolds
addressed by Casson and Gordon and the case of bounded manifolds
addressed in Theorem \ref{weakly-reducible}.  Conclusions (2) and
(3) in the above theorem do {\it not} imply that the splitting is
reducible. A Heegaard splitting defines a partition of the
boundary components of the manifold. Reducing (destabilizing) a
Heegaard splitting  does not change this partition of the boundary
components, whereas compression along $\widehat{\D_1}$
or$\widehat{\D_2}$ may change the partition. This is seen in the
following example.

\begin{ex}
Consider $M$, the exterior of the three component chain pictured
in Figure \ref{f-wr}.   (The manifold $M$ may also be though of as
$P \times S^1$, where $P$ is a pair of pants).

\begin{figure}[ht!]
{\epsfxsize = 3 in \centerline{\epsfbox{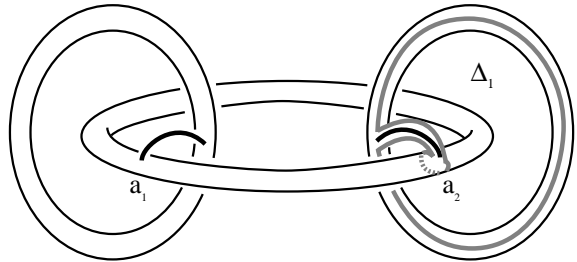}} } \caption{A
weakly reducible tunnel system for the exterior of the three
component chain.} \label{f-wr}
\end{figure}

It is not difficult to see that when a neighborhood of the arcs
$a_1$ and $a_2$ are tunneled out of $M$, a handlebody is the
result.  Thus, $\{a_1 \cup a_2 \}$ is a tunnel system for $M$.
Moreover, this system is weakly reducible: let $\D_1 = D_1$ be the
cocore of $a_1$ and $\Delta_2$ be the compressing disk for $M-
N(a_1 \cup a_2)$ whose boundary is indicated in the figure,
running over $a_2$ twice. Since, $M$ does not contain any closed
non-peripheral incompressible surfaces, Theorem
\ref{weakly-reducible} implies that this splitting can be
compressed to a splitting of lower genus. However, the tunnel
system $\{a_1 \cup a_2\}$ cannot possibly be reducible, three is
the minimal genus of a Heegaard splitting for which  all three
boundary components of $M$ are contained in the same compression
body.  In fact, either of the arcs $a_1$ or $a_2$ taken alone are
a handle system for $M$.  This induces a genus 2 Heegaard
splitting of $M$ where one compression body contains two boundary
components of $M$ and the other compression body contains one
boundary component of $M$.
\end{ex}

We now refine these methods to address the problem outlined in
Section \ref{s-main}.    The setup is the same as in that section:
$M$ is a closed manifold, $\F$ is a complete collection of
surfaces, $\G$ is the spine of a Heegaard splitting that has been
arranged to intersect $2~\F$ minimally, and $\O_i$ is the support
of a relatively outermost disk $D_i$ for some component manifold
$M_i$. The proof uses the notation and closely follows the proof
of Theorem \ref{weakly-reducible}.

\begin{prop}
\label{single-handle} If $h(\O_i)=1$ then $\O_i$ is a handle
system for $M_i$.  In particular,  $h(M_i) \leq 1$.
\end{prop}

\begin{proof}
By Theorem \ref{t-tunnelsys} we know that $\G_i$ is a tunnel
system for $M_i$; $H = \bdy N(\G_i \cup \bdy M_i)$ is a Heegaard
surface for $M_i$. We may write $M_i = C_1 \cup_H C_2$ where $C_1$
is a compression body containing {\it all} components of $\bdy
M_i$ and $C_2$ is a handlebody.

\begin{figure}[ht!]
{\epsfxsize = 4 in \centerline{\epsfbox{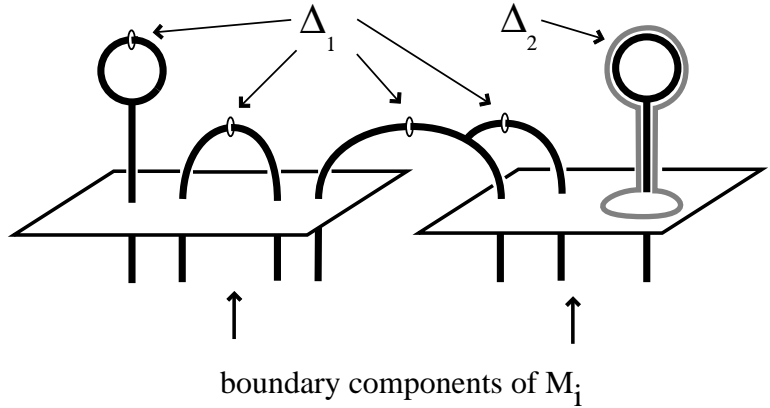}} }

\caption{The tunnel system $\G \cap M_i$ is weakly reducible.}
\label{f-d1d2}
\end{figure}

If $h(\G_i)=1$ then $\O_i = \G_i$ and the result holds trivially.
So we assume that $h(\G_i) > 1$.  This implies that $M_i = C_1
\cup_H C_2$ is weakly reducible:  let $\D_1 \subset C_1$ be a
complete collection of compressing disks for  $N(\G_i - \O_i)$,
and $\D_2 = D_i$, the relatively outermost disk.  See Figure
\ref{f-d1d2}. We may choose $\Delta_1$ so that $Y_1$ does not
contain any balls, every component is attached to $\bdy M$.

Following the proof of Theorem \ref{weakly-reducible}, by further
compressing $H_0$ we obtain a sequence of weak reducing systems,
with strictly decreasing complexity. Let  $\D_0'=\D_1' \cup \D_2'$
be the final disk system; for this system the surface $H_0'$
consists of 2-spheres and incompressible surfaces.

\begin{claim} Every incompressible component of $H_0'$ is
peripheral.
\end{claim}

Otherwise, by Theorem \ref{weakly-reducible}  $M_i$ contains an
orientable, non-peripheral incompressible surface $S$.   Recall
that $H_1$ is the surface obtained by compressing $H$ along
$\D_1$, in this case it consists of copies of boundary components
(possibly none), and either one or two boundary components with a
handle, $\bdy N(\O_i)$, attached. To obtain $S$, we further
compress along $\D_2=D_i$ and perhaps along additional disks.
Therefore $S$ has genus less than or equal to the genus of a
boundary component when $\O_i$ is attached to a single boundary
component, or strictly less than the sum of the genera of two
boundary components when $\O_i$ is attached to two boundary
components. This violates our assumption that the decomposition
along $\F$ was complete.  This completes the proof of the claim.

\bigskip

We conclude, as in Theorem \ref{weakly-reducible}, that  some
component $M' \subset X_1'$ or $M' \subset X_2'$ is a copy of
$M_i$ with some balls removed.

\begin{claim}
If $M' \subset X_1'$ then $\O_i$ is a handle system for $M_i$.
\end{claim}

In this case some component of $H_1' = \bdy Y_1'$ is a Heegaard
surface for $M_i$.  Recall that $Y_1'$ was obtained by first
compressing $C_1 =N(\bdy M_i \cup \G_i)$ along $\D_1$ yielding
$Y_1 = (\bdy M_i \times I) \cup N(\O_i) $, and then perhaps
compressing further. But, the only compressions remaining are
along the cocore of $N(\O_i)$ or the separating disk which is the
double of the cocore (only if $\O_i$ is attached to a single
component of $\bdy M_i$).  But, we could not have compressed along
either of these: compressing along the cocore leaves $Y_1' = \bdy
M_i \times I$ whose boundary cannot include a Heegaard surface
($M_i$ is not a compression body), and compressing along the
double would imply that the Heegaard surface is the boundary of a
solid torus, in particular $M_i$ has genus 1. This is not possible
since $\bdy M_i$ has positive genus. Thus $\Delta_0$ is the final
weak reducing system and $H_1' = H_1$. The Heegaard surface is the
boundary of the component of $Y_1$ that includes $N(\O_i)$. In
other words, $\O_i$ is a handle system for $M_i$.  (Recall Figure
\ref{f-weakred}).  This completes the claim.

\bigskip

The theorem will follow from the proof of the following claim. In
it we argue that in fact, the initial disk system $\D_0$ is
``almost'' the final system $\D_0'$.  There may be one additional
compression, but, it can be controlled.

\begin{claim}
If $M' \subset X_2'$ then the intersection of $\Gamma$ with $2~\F$
is not minimal (a contradiction).
\end{claim}

Since $M'$ is isotopic to $M$ (modulo balls),  for each boundary
component, $M'$ either contains that boundary component or a
parallel copy of that boundary component.  In fact, each must be a
parallel copy: the original Heegaard surface $H$ separated the
handlebody $C_2$ from the boundary components, and then so must
$H_0'$ separate $X_2'$ from the boundary components.  (Figure
\ref{f-weakred}).

Since $M' \subset X_2'$ we have that $H_0' = \bdy X_2'$ contains
at least one parallel  copy of each boundary component of $M_i$.
It follows that $$c(H_0') \geq c(\bdy M_i).$$
We now show that the surface $H_0$ contains a parallel copy of
each component of $\bdy M_i$.  Denote the component(s) to which
the handle $\O_i$ is attached by $\bdy_1M_i$, and the others by
$\bdy_2M_i$.  The surface $H_1$ was obtained by compressing along
$\D_1$ and therefore contains a copy of each component of
$\bdy_1M_i$ and a copy of $\bdy N(\O_i) \cup \bdy_2M_i$.    It
follows that
$$c(H_1) =
  \begin{cases}
    c(\bdy M_i) + 2  & \text{ if } |\bdy_1M_i| =1 , \\
    c(\bdy M_i) + 1  & \text{ if } |\bdy_1M_i| =2 .
  \end{cases}$$
The complexity of $H_0$ is less by one if $\D_2 = D_i$ separates
and less by two if $\D_2 = D_i$ does not separate.  Unless $\O_i$
is attached to a single component and $D_i$ separates, we have
$c(H_0) \leq  c(\bdy M_i)$.   But we know that $c(\bdy M) \leq
c(H_0') \leq c(H_0)$.  This implies that $H_0=H_0'$ and by
previous comments contains a copy of each boundary component.

The remaining case is that $\bdy_1 M_i$ is a single component and
$\D_2$ is separating.  In this case $H_0$ contains $\bdy_2M_i$ and
two other components, $S'$ and $S''$, each with positive genus. In
particular, $c(H_0)= c(\bdy M_i) + 1$. There can only be one
additional compression, say to $S'$. But this implies that $S''$
is incompressible and hence is a copy of $\bdy_1M_i$.  In this
case $S''$ is a compressible torus.

We now know that the surface $H_0$ is a single parallel copy of
each boundary component and possibly a single compressible torus
(see the proof of the above claim).   The disks $\Delta_1$ were
chosen so that $C_1$, and hence $X_1$, contains at most one
component for each component of $\bdy M_i$. Therefore the torus,
if it exists, is contained in a product neighborhood of the
boundary and some component of $X_2$ is a manifold $M'$ isotopic
to $M$. See Figure \ref{f-x1x2}.

\begin{figure}[ht!]
{\epsfxsize = 2 in \centerline{\epsfbox{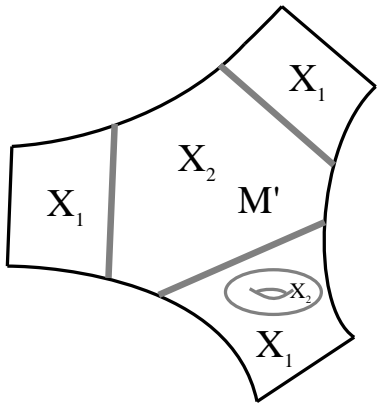}} } \caption{A
component of $X_2$ is isotopic to $M$.} \label{f-x1x2}
\end{figure}

How does the spine $\G$ intersect $M'$ ? As $M'$ is a component of
$X_2$ it is formed by attaching the 2-handles $N(\D_1)$ to a
handlebody $C_2$ which is disjoint from $\G$, again Figure
\ref{f-weakred}. But $\G$ intersects $X_2$ in a single arc for
each disk in $\D_1$ and the number of these is strictly less than
$h(\G \cap M_i)$. Then $\G$ intersects $M' \subset X_2'$ in a
subset of these arcs (perhaps all). By shrinking $M_i$ to $M'$ we
reduce $h(\G \cap M_i)$ to $|\D_1|$ which is less by at least one
(the handle $\O_i$). This contradicts the fact that $\sum_1^n h(\G
\cap M_i)$ was chosen to be minimal. Note that this may increase
the handle number of the intersection of $\G$ with the product
manifold adjacent to $M_i$. Also note that the situation in this
claim is precisely the situation that occurs when the outermost
disk demonstrates that its  support is parallel to the boundary
surface, recall Figure \ref{f-parallel}. This completes the proof
of the claim and the theorem.
\end{proof}

\begin{rem}
It is in the last claim of the proof that we are using the fact
that $h(\O_i) = 1$.  It works because we have one handle, $\O_i$,
and one compressing disk for the complement, $\D_2=D_i$.  This
implies that the original disk system $\D_0$ is ``almost'' the
final compressing system $\D_0'$. If the handle number were
greater than one then we would have a discrepancy between the
handle number of $\O_i$ and the number of compressions in $\D_2$,
we would need to compress $C_2$ further, and be forced to change
from the original disk collection $\D_1$ to a new disk collection
$\D_1'$. We cannot directly conclude that $\G$ intersects each
disk of $\D_1'$ exactly once and the intersection of $\G$ with
$M'$ might not be lower than that with $M_i$.
\end{rem}

\Addresses\recd

\end{document}